\input amstex
\documentstyle{amsppt}
\def\pf{\hfill\hfill\qed}

\def\br{\Bbb R}

\def\bz{\Bbb Z}

\leftheadtext{Augustin Banyaga}
\TagsOnRight
\NoBlackBoxes
\topmatter
\title
EXAMPLES OF NON $d_{\omega}$-EXACT LOCALLY CONFORMAL SYMPLECTIC FORMS
\endtitle
\author
Augustin Banyaga 
\endauthor

\abstract

We exhibit two three-parameter families of locally conformal symplectic forms on the solvmanifold $M_{n,k}$
considered in [1], and show, using the Hodge-de Rham theory for the Lichnerowicz cohomology that that they are not $d_{\omega}$ exact, i.e. their Lichnerowicz classes are non-trivial (Theorem 1). This has several important geometric consequences ( corollary 2, 3). This also implies that the group of automorphisms of the corresponding locally conformal symplectic structures behaves much like the group of symplectic diffeomorphisms of compact symplectic manifolds. We initiate the classification of the local conformal symplectic forms in each 3-parameter family (Theorem 2, corollary 1). We also show that the first ( and ) third Lichnerowicz cohomology classes are non-zero (Theorem 3). We observe finally that the manifolds $M_{n,k}$ carry several interesting foliations and Poisson structures.

\endabstract
\endtopmatter
\document
\baselineskip 24pt

{\bf 1. Introduction}
\vskip .1in
Using a closed 1-form $\omega$ on a smooth manifold $M$, we deform the de Rham differential $d$ on differential forms $\Lambda^*(M)$ into the operator $d_\omega$:
$$
d_{\omega}\theta = d\theta + \omega \wedge \theta ~~~~~~~~~~~~~~~\tag 1
$$
It is easy to see that $(d_\omega)^2 = 0$. Hence, we can define  $H^*_\omega (M)$  as the quotient of $Ker d_{\omega}$ by the image of $d_{\omega}$, and call it the Lichnerowicz cohomology of $M$ relative to the 1-form $\omega$ [7]. One proves that $H^*_{\omega}(M)$ is isomorphic to the cohomology $H^*(M, \Cal F_{\omega})$ of $M$ with values in the sheaf  of local functions $f$ on $M$ such that $d_{\omega}f = 0$. In [5], $H^*_{\omega}(M)$ was also characterised as the cohomology of conformally invariant forms on the minimum regular cover of $M$ over which $\omega$ pulls back to an exact form.

If the 1-form $\omega$ is exact, then the Lichnerowicz cohomology is isomorphic to the de Rham cohomology $H^*(M, \br)$ of $M$. But if $\omega$ is not exact, these two cohomologies are very different. For instance if $M$ is a compact oriented n dimensional manifold, then $H^0_{\omega}(M) \approx H^n_{\omega}(M) = 0$ for any non-exact 1-form $\omega$ [7], [8], [5].
\vskip .1in
There are few instances in which these cohomology groups are explicitely calculated. The goal of this paper is to compute some Lichnerowicz cohomology classes of the locally conformal symplectic manifold constructed by L.C.D. Andres, L.A. Cordero, M. Fernandez, and J.J. Mencia [1], herein nicknamed the "ACFM manifolds".
\vskip .1in
A {\bf locally conformal symplectic form} on a smooth manifold of dimension at least 4, is a non-degenerate 2-form $\Omega$ such that there exists a closed 1-form $\omega$ satisfying:
$$
d\Omega = - \omega\wedge\Omega ~~~~~~~~~~~~~~~~~~~~~~~~~\tag 2
$$
The closed 1-form $\omega$ above is uniquely determined by $\Omega$ and is called the Lee form of $\Omega$ [9]. The uniqueness of the Lee form is a consequence of the following elemetary fact [5]
\vskip .1in
\proclaim
{Lemma 0}
If a 2-form $\Theta$ has rank at least 4 at every point, and $\alpha$ is any 1-form, then $\alpha \wedge \Theta = 0$ implies that $\alpha$ is identically zero.
\endproclaim
\vskip .1in
If $\Omega$ is a locally conformal symplectic form with Lee form $\omega$, then equation (2) says that $d_{\omega}\Omega = 0$. Hence $\Omega$ represents an element $[\Omega] \in H^2_{\omega}(M)$.

All the examples known to the author so far were locally ( non global) conformal symplectic forms $\Omega$ with trivial Lichnerowicz class $[\Omega]$ [5], [6], [8], [13].

We exhibit  here examples of locally (non global) conformal symplectic forms with non trivial Lichnerowicz classes. 
\vskip .2in
{\bf 2. The ACFM manifold [1] and statement of the results}
\vskip .1in
Let $G_k$ be the group of matrices of the form:
$$
A =\left( \matrix
e^{kz} & 0 & 0 & x\\
0 & {e^{-kz}} & 0 & y\\
0 & 0 & 1 & z\\
0 & 0 & 0 & 1
\endmatrix \right)
$$

where $x, y, z \in \br$ and $k \in \br$ is such that $e^k + e^{-k} \in \bz - \{2\}$. The group $G_k$ is a connected solvable Lie group with coefficients: $x(A) = x$, $y(A) = y$ , $z(A) = z$ and a basis of right invariant 1-forms:
$$
dx - kx dz , ~~~~~~~~  dy + ky dz, ~~~~~~   dz 
$$
There exists a discrete subgroup $\Gamma_k \subset G_k$ such that $N_k = G_k/ \Gamma_k$ is compact [2]. The basis (3) descends to a basis of 1-forms:
$$
 \alpha, \beta, \gamma ~~~~~~~~~~~~~~~~~~~~~~~~~~~~~~~~~~~\tag 3
$$
on $N_k$. 

\vskip .1in

The forms $\gamma$ and $\alpha \wedge \beta$ are closed and their cohomology classes $[\gamma]$, $[\alpha\wedge\beta]$ generate $H^1(M, \br)$ and $H^2(M, \br)$ respectively.

Let $\lambda \in \br$  be a number such that $\lambda[\alpha\wedge\beta] \in H^2(M,\bz)$.

\vskip .1in

\proclaim
{Definition}

For any non zero $n \in \bz$, the ACFM manifolds $M_{k,n}$  are the total spaces of principal  $S^1$ bundles over $N_k$ with Chern class $n\lambda[\alpha\wedge\beta]$.

\endproclaim

\vskip .1in

For simplicity, we will denote the pull back to $M_{k,n}$  (by the projection $M_{k,n} \to N_k$) of a form $\theta$ on $N_k$ , again by $\theta$.

\vskip .1in
If $\eta$ is a connection form on $M_{k,n}$ with curvature $n\lambda (\alpha \wedge\beta)$, then:
$$ 
d\eta = n\lambda(\alpha \wedge \beta). ~~~~~~~~~~~~~~~~~~~~~~~\tag 4
$$
The set 
$$\{\alpha, \beta, \gamma, \eta\}~~~~~~~~~~~~~~~~~~~~~~~~~\tag 5$$
form a basis of 1-forms on the 4-dimensional manifold $M_{k,n}$.
\vskip .1in
Our main results are contained in theorems 1, 2 , 3 below: 
\vskip .1in
\proclaim
{Theorem 1}

(i) For each real numbers $(t_1, t_2, t_3)$ , with $t_1. t_2 \neq 0$ and $t_1t_2 \neq nk\lambda t_3^2$ the 3-parameter family of 2 -forms
$$
\Omega_{(t_1, t_2, t_3)} = t_1(\alpha \wedge \eta) + t_2(\beta \wedge \gamma) + t_3(n\lambda (\alpha\wedge \beta) - k \gamma \wedge \eta)
$$
are locally conformal symplectic forms, with the same Lee form $\omega = -k\gamma$, and their Lichnerowicz classes $[\Omega_{(t_1, t_2, t_3)}] \in H^2_{\omega}(M_{n,k})$ are non-trivial.

(ii) For each real numbers $(s_1, s_2, s_3)$ , with $s_1.s_2 \neq 0$ and $ s_1s_2 \neq n\lambda k s_3^2$, the 3-parameter family of 2 -forms
$$
\Omega_{(s_1, s_2, s_3)} = s_1(\beta \wedge \eta) + s_2(\alpha\wedge \gamma) + s_3(n\lambda (\alpha\wedge \beta) + k \gamma \wedge \eta)
$$
are locally conformal symplectic forms, with the same Lee form $\omega' = k\gamma$, and their Lichnerowicz classes  $[\Omega'_{(s_1, s_2, s_3)}] \in H^2_{(-\omega)}(M_{n,k})$ are non-trivial.             
\endproclaim
\vskip .1in
{\bf Remark 1}
\vskip .1in
The local conformal symplectic forms above are not global conformal symplectic since their Lee form is not exact [5]. Besides, $M_{n,k}$ can not carry a symplectic form since $H^2(M_{n,k},\br) = 0$ [1].
\vskip .1in
{\bf Remark 2}
\vskip .1in
The 2-forms 
$$
\Omega =n\lambda (\alpha\wedge \beta) - k \gamma \wedge \eta = d_{\omega}(\eta)
$$
$$
\Omega' =n\lambda (\alpha\wedge \beta) + k \gamma \wedge \eta= d_{\omega'}(\eta)
$$
are exact locally conformal symplectic forms. They are obviously non-degenerate
and $d_{\omega}\Omega = (d_{\omega})^2\eta = 0$. Same thing with $\Omega'$. Since the forms in theorem 1 are not exact, we have the following
\vskip .1in
\proclaim
{Corollary 1}

The locally conformal symplectic form $\Omega$ , resp $\Omega'$ is not equivalent to $\Omega_{(t_1, t_2, t_3)}$, resp. $\Omega_{(s_1, s_2, s_3)}$, i.e. there are no diffeomorphisms $\phi$ of $M_{n,k}$ and smooth function $f$ such that $\phi^*\Omega = f \Omega_{(t_1, t_2, t_3)}$ ( similarly for $\Omega'$ and
$\Omega_{(s_1, s_2, s_3)}$),

\endproclaim
\vskip .1in
\proclaim
{Corollary 2}

No Lie group can act transitively on $M_{n,k}$ preserving $\Omega_{(t-1, t_2, t_3)}$ or $\Omega_{(s_1, s_2, s_3)}.$
\endproclaim
\vskip .1in
\proclaim
{Theorem 2}

If $(t_1, t_2,t_3)\in \br^3$  satisfy  $t_2 = e^u$ and $ t_1.e^{ut} \neq nk\lambda st_3^2$ for all $t,s \in [0,1]$, there exists a  family of diffeomorphisms $\phi_t$ with $\phi_0 = id$ and
$$
\phi_t^*(\Omega_{(t_1, t_2, t_3)}) = f_t ((t_1(\alpha\wedge\eta) + \beta \wedge \gamma).
$$
We have a similar statement for $\Omega_{(s_1,s_2, s_3)}$.
\endproclaim
\vskip .1in
{\bf Question 1}
\vskip .1in
 Consider the open subset $U\subset \br^3$ defined by $ U = \{(t_1, t_2,t_3) | t_1t_2 \neq nk\lambda t_3^2\}$. The hypothesis in theorem 2 means that the point $p =( t_1,t_2,t_3) \in \br^3$ and the point $q = (t_1,1,0)$
belong to the same path component of $U$.

If $p, q\in U$ belong to two different path components, are  $\Omega_{(p)}$ and $\Omega_{(q)}$ still equivalent?

Here we can not use Moser's theorem, because we do not have a path of locally conformal symplectic forms. The path degenerates when we go from one component of $U$ to another.

\vskip .1in

For instance , take $t_1 = n\lambda/k$, $t_2 = 1$, then $t_3$ must be different from $-1/k$ and $+1/k$.
All the locall conformal symplectic forms $\Omega_{(n\lambda/k, 1, t_3)}$ are equivalent for $t_3$ in the interval $[ (-1/k) + \epsilon, (1/k) - \epsilon]$, and $\epsilon$ a small positive number.

Are the forms $\Omega_{(n\lambda/k,1,0)}$ and $\Omega_{(n\lambda/k, 1, (k+1)/k)}$ equivalent?
\vskip .1in
{\bf Remark 3}
\vskip .1in
The locally conformal symplectic form $\Omega = (n\lambda/k)(\alpha\wedge \eta) + \beta\wedge \gamma$ was found in [1] where the authors showed that it is a locally conformal Kaehler form, with non parallel Lee form, with respect to the metric $g = \alpha^2 + \beta^2 + \gamma^2 + \eta^2$.

In [10], it is proved that if a compact manifold $M$ has a riemannian metric with respect to which a nowhere zero closed 1-form $\omega$ is parallel, then $H^*_{\omega}(M) = 0$. This result and theorem 1 imply a stronger result:

\vskip .1in

\proclaim
{Corollary 3}

There is no riemanian metric on $M_{n,k}$ with respect to which  $\pm k\gamma$ is parallel.

\endproclaim

\vskip .1in

Vaisman and Goldberg [15] have found sufficient conditions for the Lee form of a compact locally Kaehker manifold to be parallel:

{\bf Theorem}
\vskip .1in
{\it Let $M, J, g)$ be a compact locally Kaehler manifold with Lee form $\omega$ a nowhere vanishing form.
If the Ricci tensor of $(M, g)$ is positive semi-definite and vanishes only in the direction of $B$, where $B$ is the vector field defined by $i(B)\Omega = \omega$, and $\Omega$ is the locally conformal Kaehler form, then $\omega$ is parallel}.

\vskip .1in
{\bf Question 2}

Let $\Omega$ be the locally conformal Kaehler form in the Goldman-Vaisman theorem above, then $H^*_{\omega}(M) = 0$. In particular $\Omega$ is $d_{\omega}$-exact.

Is Goldberg-Vaisman theorem above still true for locally conformal symplectic manifolds? 

This would give a sufficient condition for the vanishing of the Lichnerowicz cohomology and in particular the $d_{\omega}$-exactness of the locally conformal symplectic form.

\vskip .1in

\proclaim
{Theorem 3}

 The forms $\alpha$, $\alpha\wedge \eta$ and $\alpha \wedge \gamma \wedge \eta$ represent non zero elements in $H^i_{\omega}(M_{n,k})$, $ i = 1, 2, 3$.

The forms $\beta$, $\beta \wedge \eta$, $\beta\wedge\gamma\wedge\eta$ represent non zero elements in $H^i_{(-\omega)}(M_{n,k})$, $ i = 1, 2, 3$.

The classes $[\Omega_{(t_1,t_2,t_3)}]$ in theorem 1 coincide with $t_1[\alpha\wedge\eta]$. Similarly, $[\Omega_{(s_1, s_2, s_3)}] = s_1[\beta\wedge \eta]$.

\endproclaim

\vskip .1in
{\bf Question 3}
\vskip .1in
Theorem 3 says that the dimension of $H^i_{\{\omega, -\omega\}}(M_{n,k})$, for i = 1,2,3 is at least 1. Is it exactly one ?

\vskip .1in

{\bf Remarks on the group of automorphisms and the Lie algebra of infinitesimal automorphisms}

\vskip .1in

Two locally conformal symplectic forms $\Omega$ and $\Omega'$ are said to be equivalent if there exists a no-where zero function $f$ such that $\Omega' = f\Omega$. 

 A locally conformal symplectic structure $\Cal S$ is an equivalence class of locally conformal symplectic forms. Let $\Cal S_t$, respectively $\Cal S_s$ be the locally conformal symplectic structures on $M_{n,k}$ represented by the forms $\Omega_{(t_1, t_2, t_3)}$ respectively $\Omega_{(s_1, s_2, s_3)}$, and let $G(\Cal S_t)$, $G(\Cal S_s)$ be the corresponding automorphism groups, i.e. the group of all diffeomorphisms $\phi$ of $M_{n,k}$ such that $\phi^*(\Omega_{(t_1, t_2, t_3)}) = f\Omega_{(t_1, t_2, t_3)}$ for some no-where zero function $f$ (same definition for $\Omega_{(s_1, s_2, s_3)}$).
\vskip .1in
Let $\Cal L(\Cal S_{t}(M_{n,k}))$ and $\Cal L(\Cal S_{s}(M_{n,k}))$  be the Lie algebra of infinitesimal automorphisms of the locally conformal symplectic structures $\Cal S_t$ and $\Cal S_s$. This is the Lie algebra of vector fields $X$ on $M_{n,k}$ such that $L_X\Omega_{(t_1,t_2,t_3)} = f \Omega_{(t_1,t_2,t_3)}$, respectively  $L_X\Omega_{(s_1,s_2,s_3)} = f \Omega_{(s_1,s_2,s_3)}$,for some function $f$,
where $L_X$ stands for the Lie derivative in the direction $X$.

\vskip .1in

Let $M$ be a compact manifold equipped with a locally conformal symplectic form $\Omega$, with Lee form $\omega$  and such that $[\Omega]$ is a non trivial element of $H^2_{\omega}(M)$ , then the groups of automorphisms of the corresponding locally conformal symplectic structure, and the corresponding Lie algebra, behave very much like in the symplectic case [8]. For instance there is a Calabi invariant on the identity component of the group with values in a quotient of $H^1_{\omega}(M)$ and its kernel is a simple group [8], exactly like in the symplectic case [3]. At the level of the Lie algebras, the mapping $X\mapsto i(X)\Omega$ is a surjective Lie algebra homomorphism into $H^1_{\omega}(M)$ whose kernel is the commutator subagebra [8]. The fact that $i(X)\Omega$ is $d_{\omega}$ closed comes from the following fact:

\proclaim
{Corollary 4}

For every vector field $X \in \Cal L(\Cal S_{t}(M_{n,k}))$ , then 
$$
L_X\Omega_{(t_1,t_2,t_3)} = k\gamma(X)\Omega_{(t_1,t_2,t_3)}.
$$
For $X \in \Cal L(\Cal S_{s}(M_{n,k}))$ , then 
$$
L_X\Omega_{(s_1,s_2,s_3)} = - k\gamma(X)\Omega_{(s_1,2_2,s_3)}.
$$
\endproclaim

\vskip .1in

An immediate consequence of corollary 4 is the following fact: for any infinitesimal automorphism $X \in
\Cal L(\Cal S_{t}(M_{n,k}))$ or in  $X \in \Cal L(\Cal S_{s}(M_{n,k}))$ :
$$
\int_{M_{n,k}} \gamma(X)(\Omega_{(t_1,t_2,t_3)})^2  = \int_{M_{n,k}} \gamma(X)(\Omega_{(s_1,s_2,s_3)})^2 = 0.
$$

\vskip .2in

{\bf 3. Proofs of the results}

\vskip .1in

The proofs rest on the Hodge-de Rham theory for the $d_{\omega}$ operator and lemma 1 below.

{\bf The Hodge-de Rham theory for the $d_{\omega}$ operator [7]}
\vskip .1in
Let $M,g)$ be a compact oriented n-dimensional riemannian manifold and $\omega$ a 1-form on $M$. Let $*$ be the Hodge-de Rham star operator defined by $g$ and $\delta$  the codifferential:
$ \delta (\theta) = (-1)^{(nl + n +1)}(* \circ d \circ *)\theta$ for $\theta \in \Lambda^l(M)$, and  let $\Cal U_{\omega}$ be the operator: $\Cal U_{\omega}(\theta )= (-1)^{nl + n}(* (\omega\wedge * \theta))$ for $\theta \in \Lambda^l(M)$.
\vskip .1in
We let
$$
\delta_{\omega} = \delta + \Cal U_{\omega}
$$
Consider the the inner product $<,.>$ on $\Lambda^l(M)$:
$$
<\rho,\nu> = \int_M \rho\wedge *\nu.
$$
One has that $<d_{\omega}\rho,\nu> = <\rho, \delta_{\omega}\nu>$. Since $M$ is compact and $(\Lambda^*(M), d_{\omega})$ is elliptic, we get the following result, proved in [7]:
\vskip .1in
\proclaim
{Theorem}

We have the orthogonal  decomposition:
$$
\Lambda^l(M) = (\Cal H^l_{\omega}(M)) \oplus (d_{\omega}(\Lambda^{l-1}(M))) \oplus (\delta_{\omega}(\Lambda^{l+1}(M)))
$$
where
$$
\Cal H^l_{\omega}(M) = \{ \theta \in \Lambda^l |d_{\omega} \theta = \delta_{\omega}\theta = 0 \}
$$
is the space of $\omega$-harmonic forms.

Fianally, we have:
$$
H^l_{\omega}(M) \approx \Cal H^l_{\omega}(M). 
$$
\endproclaim
\vskip .1in
On the manifold $M_{n,k}$, we will put the riemannian metric $g$ which makes the forms $\alpha, \beta,\gamma, \eta$ orthonormal, i.e. $g = \alpha^2 + \beta^2 + \gamma^2 + \eta^2$. We will consider the closed 1-forms $\omega = -k\gamma$ and $(-\omega)$ respectively.
\vskip .1in
\proclaim
{Lemma 1}

If $\alpha, \beta, \gamma, \eta$ is the basis (5) of 1-forms on $M_{k,n}$, we have:
$$
d\alpha = - k \alpha \wedge \gamma ~~~~~~~~~~~~~~~~~~~~~\tag 6
$$
$$
d\beta = k \beta \wedge \gamma~~~~~~~~~~~~~~~~~~~~~~~~~~~\tag 7
$$
\endproclaim

\vskip .1in

\demo
{Proof}

Let $\{X, Y, Z, T\}$ be the basis of global vector fields on $M_{k,n}$ dual to the basis of 1-form (5).
A general 2-form $\theta$ on $M_{k,n}$ is uniquely written as : 

$\theta = A(\alpha \wedge \beta) + B (\alpha\wedge \gamma )+ C(\alpha \wedge \eta )+ D (\beta \wedge \gamma) + E (\beta \wedge \eta) + F (\gamma \wedge \eta)$,

 where

$$A = \theta(X,Y),~~~~~~~~~~~~~ B = \theta(X, Z),$$
$$ C = \theta(X, T),~~~~~~~~~~~ D = \theta(Y, Z),$$
$$ E = \theta(Y, T),~~~~~~~~~~~ F = \theta(Z, T).$$

If $\rho$ denotes any of the basic 1-forms $\alpha, \beta, \gamma, \eta$ and $\xi, \xi'$ any pair
of the vector fields $\{X, Y, Z, T\}$ dual to the basis of 1-forms above, then:
$$
(d\rho)(\xi, \xi') = - \rho([\xi,\xi']) + \xi.\rho(\xi') -\xi'.\rho(\xi)  = - \rho([\xi,\xi']) ~~~~~~~~~~~~~~~~~~~~~~~~~~~~~\tag 8
$$
since $\rho(\xi), \rho(\xi')$ are 1 or 0, and hence their directional derivatives are zero.

\vskip .1in

We have the following facts [1]:
$$
[X, Z] = k X
$$
$$
[X, Y] = -n\lambda T
$$
$$
[Y, Z] = - kY
$$
$$
[X, T] = [Y, T] = [Z, T] = 0 
$$
Applying the above facts to $\theta = d\alpha$, and using the commutation relations above, we see that the only non zero coefficient is $B = - k$. Hence
$$
d\alpha = -k \alpha \wedge \gamma
$$
A similar inspection for $d\beta$ yields
$$
d\beta = k\beta \wedge \gamma
$$,
\pf
\enddemo

\vskip .1in

\proclaim
{Proposition 1}

1. Let $\omega = - k\gamma$, then  $\alpha$, $\alpha \wedge \eta$, $\alpha \wedge \gamma \wedge \eta$ are $\omega$-harmonic.

The form $\beta \wedge\gamma$ is $d_{\omega}$-exact.  In fact, $d_{\omega}((1/2k) \beta) = \beta\wedge \gamma$.

2. Let $\omega' = k\gamma$, then $\beta$, $\beta \wedge \eta$, $\beta \wedge \gamma \wedge \eta$
are $\omega'$-harmonic.

 The form $\alpha \wedge \gamma$ is $d_{\omega'}$-exact. In fact, $d_{(-\omega)}((-1/2k) \alpha) = \alpha\wedge \gamma$.

\endproclaim

\vskip .1in

\demo
{Proof}

By lemma 1, $ d\alpha = -k\alpha \wedge \gamma = - \omega \wedge \alpha$. Hence $d_{\omega}\alpha = 0$. Using (4), we get: 

$d(\alpha \wedge \eta) = - k\alpha \wedge\gamma\wedge\eta - \alpha\wedge (n\lambda \alpha \wedge\beta)=  - k\alpha \wedge\gamma\wedge\eta = -\omega \wedge (\alpha \wedge\eta)$. Therefore $d_{\omega}(\alpha\wedge\eta) = 0$.

$d(\alpha\wedge\gamma\wedge\eta) = -k\alpha\wedge\gamma\wedge\gamma\wedge\eta - \alpha\wedge\gamma \wedge(n\lambda \alpha\wedge \beta) = 0$ and $\omega \wedge (\alpha\wedge\gamma\wedge\eta) = 0.$ Hence $d_{\omega}(\alpha\wedge\gamma\wedge \eta) = 0$.

We just have proved that the forms listed in (1) are $d_{\omega}$ closed. Let us now show that $\delta_{\omega}$ of these forms are zero:

Up to a sign, $*\alpha = \beta\wedge\gamma\wedge\eta$, hence $d*\alpha = k\beta \wedge\gamma\wedge \gamma \wedge\eta + \beta\wedge \gamma\wedge (n\lambda \alpha \wedge \beta) = 0$. Hence $\delta(\alpha) = 0$. On the other hands, $\omega \wedge*\alpha = - k\gamma\wedge \beta\wedge\gamma\wedge\eta = 0$, i.e. $\Cal U_{\omega}(\alpha) = 0$, and hence $\delta_{\omega}\alpha = 0$.

Up to a sign, $*(\alpha \wedge \eta) = \beta \wedge \gamma$, so $d* (\alpha \wedge \gamma) =  k\beta \wedge \gamma \wedge \gamma = 0$. This Shows that $\delta (\alpha\wedge\eta) = 0$.

We have that $\omega *(\alpha \wedge\eta) = - k\gamma \wedge\beta\wedge\gamma =0$. Hence $\Cal U_{\omega} = 0$, and thus $\delta_{\omega}(\alpha\wedge\eta) = 0$.

Up to a sign, $*(\alpha \wedge \gamma \wedge \eta) = \beta$. Hence $* d *(\alpha \wedge \gamma \wedge \eta) = *(k\beta\wedge \gamma) = k \alpha \wedge \eta$. Hence $\delta(\alpha\wedge \gamma\wedge \eta) = (-1)^{17}k\alpha \wedge\eta = - k\alpha\wedge\eta$.

 On the other hands $\Cal U_{\omega}(\alpha\wedge \gamma\wedge \eta) =  +* ((-k\gamma)\wedge *(\alpha \wedge \gamma \wedge \eta)) = *(-k \gamma\wedge \beta) = *(k\beta\wedge \gamma) = k\alpha\wedge\eta$. Hence $\delta_{\omega}(\alpha\wedge\gamma\wedge\eta) = 0$.

Finaly, $d_{\omega}\beta = k\beta \wedge \gamma + (-k\gamma)\wedge\beta = 2k(\beta\wedge\gamma)$.
Therefore $ \beta\wedge\gamma = d_{\omega}((1/2k) \beta)$.

This proves (1). Similar calculations using $\omega' = + k\gamma$ yields (2). \pf

\enddemo
\vskip .1in
{\bf End of the proof of theorem 1}
\vskip .1in
Observe first that 
$$
n\lambda \alpha \wedge \beta - k\gamma \wedge \eta = d\eta + \omega \wedge \eta = d_{\omega}\eta
$$
($\omega = -k\gamma$). Therefore:
$$
\Omega_{(t_1,t_2,t_3)} = t_1 (\alpha\wedge\eta) + d_{\omega}((t_2/2k) \beta + t_3 \eta) ~~~~~~~~\tag 9
$$
Since $d_{\omega}(\alpha \wedge\eta) = 0$ and $(d_{\omega})^2 = 0$, it follows that $\Omega_{(t_1,t_2,t_3)}$ is $d_{\omega}$-closed.

Similarly,
$$
\Omega_{(s_1,s_2,s_3)} = s_1(\beta\wedge\eta) + d_{\omega'}((-s_2/2k) \alpha + s_3\eta). ~~~~~~~~~~~~~~~~~~~~~~~~~\tag 10
$$
hence $\Omega_{(s_1,s_2,s_3)}$ is $d_{\omega'}$-closed ( with $\omega' = k\gamma$).

The Lichnerowicz cohomology classes of these forms $[\Omega_{(t_1,t_2, t_3)}]$ are  $t_1[\alpha\wedge \eta]$ and
$[\Omega_{(s_1,s_2, s_3)}] = s_1[\beta\wedge \eta]$, which are non-trivial, since $\alpha\wedge\eta$ and $\beta\wedge\eta$ are non-zero harmonic forms by proposition 2. By the Hodge-de Rham theory, they represent non trivial elements of $H^2_{\omega}(M_{n,k})$, respectively $H^2_{(-\omega)}(M_{n,k}).$ 

\vskip .1in

The forms in (i) and (ii) are non-degenerate. Indeed an immediate calculation gives:

$(\Omega_{(t_1,t_2,t_3)})^2 = 2(t_1t_2 -nk\lambda t_3^2)(\alpha \wedge \beta \wedge \gamma \wedge \eta) \neq 0$ iff $t_1t_2 -nk\lambda t_3^2 \neq 0$

$(\Omega_{(s_1,s_2,s_3)})^2 = -2(s_1s_2 -nk\lambda s_3^2)(\alpha \wedge \beta \wedge \gamma \wedge \eta)\neq 0$ iff
$s_1s_2 -nk\lambda s_3^2 \neq 0$.
This completes the proof of theorem 1. \pf 

\vskip .1in

{\bf Proof of theorem 3}

\vskip .1in

Proposition 1 shows that the forms considered in theorem 3 are $\omega$ or $\omega'$ harmonic. Hence they represent non-zero elements in the corresponding Lichnerowicz cohomologies. \pf

\vskip .1in

{\bf Proof of theorem 2}

\vskip .1in
We will use the locally conformal symplectic forms version of Moser theorem [11]  proved in [5]:

\vskip .1in
\proclaim
{Theorem}

Let $\Omega_t$ be a smooth family of locally conformal symplectic forms on a compact manifold $M$. Suppose that for all t, the Lee form of $\Omega_t$ is the same 1-form $\omega$ and $\Omega_t -\Omega_0$ is $d_{\omega}$ exact, then there exists a smooth family of diffeomorphisms $\phi_t$, with $\phi_0 = id$ and a smooth family of functions $f_t$ such that $\phi_t^*\omega_t = f_t\Omega_0$

\endproclaim

\vskip .1in

Let $(t_1,t_2,t_3)\in \br^3$ such that $t_2 = e^u$ and $t_1.e^{ut} \neq nk\lambda st_3$ for all $t,s \in [0,1]$, then
$$
 \Omega^{(s,t)}_{(t_1,t_3)} = t_1(\alpha\wedge\eta) + e^{ut}(\beta\wedge\gamma + e^{-ut}st_3 d_{\omega}\eta) = t_1\alpha\wedge\eta + d_{\omega}(e^{ut}((1/2k)\beta) + e^{-ut} st_3 \eta))
$$
is a smooth 2-parameter family of locally conformal symplectic forms, all with the same cohomology class $t_1[\alpha\wedge\eta]$, with the same Lee form $\omega = -k\gamma$.

By Moser theorem for locally conformal symplectic forms, for fixed s, there exists a smooth family of diffeomorphisms $\psi^t$ depending smoothly on the parameter $s$, such that $\psi^0 = id$ and $(\psi^t)^*\Omega^{(s,t)}_{(t_1,t_3)} = f^{(t,s)} \Omega^{(0,s)}_{(t_1,t_3)}$, for some functions $f^{(t,s)}$, depending smoothly on s and t.

Apply again Moser theorem to the family $\Omega_s =\Omega^{(0,s)}_{(t_1,t_3)}$. There is a family of
diffeomorphisms $\psi'_s$ such that $(\psi'_s)^*) \Omega_s = g_s\Omega_0 = g_s((t_1 \alpha\wedge\eta) + \beta\wedge \gamma)$.

The required family of diffeopmrphism is $\psi^t \circ \psi'_s$ : it pulls $\Omega^{(s,t)}_{(t_1,t_3)}$  back to a multiple by a function of $\Omega^{(0,0)}_{(t_1,t_3)} = t_1(\alpha \wedge \eta + \beta\wedge \gamma$.
\pf
\vskip .1in
{\bf Proof of corollary 3}

Since $\Omega_{(t_1,t_2,t_3)}$ resp. $\Omega_{(s_1,s_2,s_3)}$  are non global and not $\omega$-exact, rsp. $\omega'$-exact, proposition 2.4 of [13] implies that $M_{n,k}$, equipped with the forms above, is not a homogeneous locally conformal symplectic manifold.
\pf
\vskip .1in

{\bf Proof of corollary 4}

We need to recall the generalized Lee homomorphism [13], [6].

Let $\Theta$ be a 2-form of rank everywhere larger or equal to 4, such that there exists a closed 1-form $\eta$ safisfying
$$
d\Theta = - \eta \wedge \Theta
$$
on a smooth manifold $M$. We know that $\eta$ is uniquely determined by $\Theta$ ( lemma 0).

\vskip .1in

Let $\Cal L(\Theta, M)$ denote the Lie algebra of vector fields $X$ on $M$ such that $L_X\Theta = \mu_X\Theta$.
\vskip .1in

For $X\in \Cal L(\Theta, M)$, set $\theta = i(X)\Theta$, where i(.) is the interior product operator.

We have

$L_X\Theta = \mu_X \Theta = di(X)\Theta + i(X)d\Theta$

$= d\theta + i(X)(-\eta\wedge\Theta)$

$= d\theta - \eta(X)\Theta + \eta \wedge i(X)\Theta$

$= d\theta + \eta\wedge \theta -\eta(X)\Theta  = d_{\eta}\theta -\eta(X)\Theta$

Hence:

$$
d_{\eta}\theta = l(X)\Theta ~~~~~~~~~~~~~~~~~~~~~~~\tag 11
$$
where
$$
l(X) = \mu_X + \eta(X).~~~~~~~~~~~~~~~~~~~~~~~~~~~\tag 12
$$

\vskip .1in

\proclaim
{Proposition 2}

If the manifold $M$ is connected, then the function $l(X)$ is a constant for all $X\in \Cal L(\Theta, M)$, and the map $ l : \Cal L(\Theta, M) \to \br, ~~~~~~~~~ X\mapsto l(X)$, is a Lie algebra homomorphism

\endproclaim

\vskip .1in
A simple proof that $l(X)$ is constant goes as follows:

$0 = (d_{\eta})^2\theta = d_{\eta}(l(X)\Theta) = d(l(X)\Theta) + \eta\wedge (l(X)\Theta)$

$ = d(l(X))\wedge \Theta + l(X)(-\eta\wedge \Theta) +\eta\wedge (l(X)\Theta)$

$ = d(l(X))\wedge \Theta $

\vskip .1in

By lemma 0, $d(l(X)) = 0$. Therefore $l(X)$ is a constant, since $M$ is connected.

\vskip .1in

This homomorphism was constructed by Vaisman [13] in the case where $\theta$ is a locally conformal symplectic form.

As an immediate consequence of proposition 2 we get the following proposition, which belongs to the folklore:

\proclaim
{Proposition 3}

A locally conformal symplectic form $\Omega$ with Lee form $\omega$ on a connected manifold is $d_{\omega}$-exact iff there exists an infinitesimal automorphism $X$ such that $l(X) \neq 0$.

\endproclaim

\vskip .1in

For the convenience of the reader, we include a proof.

\demo
{Proof}

If $X$ is an infinitesimal automorphism such that $l(X) \neq 0$, set :  $\theta = (1/l(X))i(X)\Omega$, then (11) yields : $\Omega = d_{\omega}\theta$.

Converseky, if $\Omega = d_{\omega}\theta$ for some 1-form $\theta$. Define a vector field $X$ by the equation $ i(X)\Omega = \theta$. An easy calculation shows that $L_X\Omega = (1- \omega(X))\Omega$. This proves that $X$ is an infinitesimal automorphism: $L_X\Omega = \mu_X\Omega$ with $\mu_X = 1- \omega(X)$, that is $l(X) = 1$.

\enddemo

\pf

\vskip .1in

The proof of corollary 2 is an immediate consequence of theorem 1. The Lee homomorphism corresponding to our locally conformal symplectic forms is trivial. Hence  if a vector field $X$ is such that $L_X\Omega_{(t_1,t_2,t_3)} = f \Omega_{(t_1,t_2,t_3)}$, then (12) yields : $f = - \omega(X)$.
\pf

\vskip .2in

{\bf 4. Some foliations and Poisson structures on $M_{n,k}$} 

\vskip .1in

The manifold $M_{n,k}$ carries several remarkable foliations:

The codimension 1 foliations $\Cal F_1$, $\Cal F_2$, $\Cal F_3$, defined by the integrable forms $\gamma$, $\alpha$ and $\beta$.

The five 2-dimensional foliations $\Cal F_4$,... $\Cal F_8$, tangent to the involutive distributions 
$$
\{X, Y\}, \{Y,Z\}, \{X,T\}, \{Y, T\}, \{Z, T\}
$$
 
For instance the tangent space to the leaves of the foliation $\Cal F_1$ defined by $\gamma = 0$, is spanned by $X, Y, T$, The restriction of the locally conformal symplectic form $\Omega = \alpha\wedge\eta + \beta\wedge\gamma$ is the 2-form $\alpha\wedge\eta$,
whose kernel is spanned by $Y$. The quotient of each leaf by the trajectories of $Y$ is spanned by $X, T$, which define the 2-dimensional foliation  $\Cal F_7$. The 2-form $\alpha \wedge \eta$ is a symplectic form on each leaf. We thus eee that $\Cal F_7$ is a symplectic foliation.
\vskip .1in
The restriction of $\eta$ to the leaves of the foliation $\Cal F _1$ is a contact form since $\eta\wedge\eta = n\lambda \alpha \wedge \beta \wedge \eta$. This is why we have the locally conformal symplectic forms $\Omega$ and $\Omega'$ in remark 2.
\vskip .1in
On each leaf of the 2 dimensional foliations there are is an obvious symplectic forms. Hence they are all symplectic foliations.

According to Vaisman [14], a symplectic foliation gives raise to a Poisson structure. The bracket of two functions is the so-called Dirac bracket: let $\Cal F$ be a symplectic foliation on a manifold $M$, the Dirac bracket $\{f,g\}$ of two functions $f,g$ on $M$ is $\{f,g\}(x) = (\Omega_F)(x) (X_{(f_{|F})}, X_{(g_{|F})})(x)$, where $\Omega_F$ is the symplectic form on the leaf $F$ through $x$ and $(X_{(f_{|F})}$ is the sympectic gradient of the restriction $f_|F$ of $f$ to the leaf $F$, with respect to the symplectic form $\Omega_F$.

This way, we get 5 Poisson structures on $M_{n,k}$.
\vskip .1in
{\bf Question 4}
\vskip .1in
Are there any compatibility relations among these Poisson structures that can be deducted from the commutation relations of the vector fields $\{X, Y, Z, T\}$?
\vskip .1in
{\bf Remark 4}
\vskip .1in
Since $\gamma$ is a closed 1-form without zero, a theorem of Tischler[12] asserts that $M_{n,k}$ is fibered over $S^1$.

\newpage

\centerline{\bf REFERENCES}

\item{[1]} L.C. Andres, L.A. Cardero, M. Fernandez and J.J. Mencia, {\it Examples of four-dimensional compact locally conformal Kaehker solvmanifolds}, Geom. Dedicata 29(1989) 227-- 232.

\item{[2]} L. Auslander, L. Green, and F. Hahn, {\it Flows on homogeneous spaces}, Ann. of Math. Studies 53, Princeton Univ. Press, 1963.

\item{[3]} A. Banyaga, {\it Sur la structure du groupe des diffeomorphismes qui preservent une forme symplectique}, Comment. Math. Helv. 53(1978) 174-227

\item{[4]} A. Banyaga, {\it The structure of classical diffeomorphism Groups}. Kluwer Academic Publishers, Mathematics and its applications, Vol 400 (1997)

\item{[5]} A. Banyaga, {\it Some properties of locally conformal symplectic structures}, Comment. Math. Helveticii, 77(2002) 383-398.

\item{[6]} A. Banyaga, {\it On the geometry of locally conformal symplectic manifolds}, Infinite dimensional Lie Groups in Geometry and Representation Theory, A. Banyaga, J. Leslie and T. Robart , Eds. World Scientific (2002) pp 79- 91

\item{[7]} F. Guedira and A. Lichnerowicz, {\it Geometrie des algebres de Lie locales de Kirillov}, J.Math. Pures et Appl. 63(1984), 407-484.

\item{[8]} S. Haller and T. Rybickii, {\it On the group of diffeomorphisms preserving a locally conformal symplectic structure}, Ann. Global Anal. and Geom. 17 (1999) 475-502 

\item{[9]} H.C. Lee,  {\it A kind of even-dimensional differential geometry and its application to exterior calculus}. Amer. J. Math. 65(1943) pp 433-438.

\item{[10]} M. de Leon, B. Lopez, J.C. Marrero, E. Padron, {\it On the computation of the Lichnerowicz-Jacobi cohomology}, Jour. Geom. and Physics 44(2003) 507-522

\item{[11]} J. Moser, {\it On the volume element of a manifold}, Trans. AMS 120(1965) 286-294

\item{[12]} D. Tischler, {\it On fibering certain manifolds over $S^1$}. Topology 9(1970), 153-154.

\item{[13]} I. Vaisman, {\it Locally conformal symplectic manifolds}, Inter. J. Math.Sci 8(3) (1985)
521-536

\item{[14]} I. Vaisman, {\it Lectures on the geometry of Poisson manifolds}, Progress in Math, Vol 118, Birkhauser, 1994

\item{[15]} I Vaisman and S. I. Goldberg, {\it On compact locally conformal Kaehler manifolds with nonnegative sectional curvature}, Ann. Fac. SCi. Toulouse, 2(198), 117-123.

 \vskip .2in
\baselineskip 12pt
\noindent Department of Mathematics \newline
\noindent The Pennsylvania State University \newline
\noindent University Park, PA 16802

\enddocument